\documentclass[]{amsart}
\usepackage{amssymb}
\usepackage{amsmath}
\author{Raf Cluckers$^*$}
\thanks{$^*$ Research Assistant of the Fund for Scientific Research --
 Flanders (Belgium)(F.W.O.)}
\address{Katholieke Universiteit Leuven, Department of Mathematics,
Celestijnenlaan 200B, 3001 Leu\-ven, Bel\-gium }
\email{raf.cluckers@wis.kuleuven.ac.be}
\urladdr{http://www.wis.kuleuven.ac.be/wis/algebra/Raf/index.html}
\date{}
\theoremstyle{plain}
\newtheorem{theorem}{Theorem}[section]

\newtheorem{lemma}[theorem]{Lemma}

\theoremstyle{definition}
\newtheorem{definition}[theorem]{Definition}

\newtheorem{remark}[theorem]{Remark}

\newcommand{\Z}{\mathbb{Z}}

\newcommand{\Q}{\mathbb{Q}}

\newcommand{\PP}{\mathcal{P}}

\newcommand{\SSS}{\mathcal{S}}
\newcommand{\aaa}{\alpha}
\newcommand{\bbb}{\beta}

\newcommand{\ddd}{\delta}
\newcommand{\lP}{\mu P_n}

\newcommand{\Ldan}{\mathcal{L}_{\mathrm an}}

\DeclareMathOperator{\sq}{\square}

\newcommand{\OOL}{{\mathcal C}}
\newcommand{\OOLp}{{\mathcal C}_{\rm simple}}
\newcommand{\cP}{\mathcal{P}}
\newcommand{\abs}[1]{\lvert#1\rvert}

\title{Analytic p-adic Cell Decomposition and
Integrals}


\begin{document}
\begin{abstract}
We prove a conjecture of Denef on parameterized $p$-adic analytic
integrals using an analytic cell decomposition theorem, which we
also prove in this paper. This cell decomposition theorem
describes piecewise the valuation of analytic functions (and more
generally of subanalytic functions), the pieces being
geometrically simple sets, called cells. We also classify
subanalytic sets up to subanalytic bijection.
\end{abstract}
\maketitle
\section{Introduction}
Let $p$ denote a fixed prime number, $\Z_p$ the ring of $p$-adic
integers, $\Q_p$ the field of $p$-adic numbers, let $|\cdot|$
denote the $p$-adic norm and $v(\cdot)$ the $p$-adic valuation.
\par
Let $f=(f_1,\ldots,f_r)$ be an $r$-tuple of restricted power
series over $\Z_p$ in the variables
$(\lambda,x)=(\lambda_1,\ldots,\lambda_s,x_1,\ldots, x_m)$, i.e.,
the $f_i$ are power series converging on $\Z_p^{s+m}$. To $f$ we
associate a parametrized $p$-adic integral
\begin{equation}\label{integral}
I(\lambda)=\int\limits_{\Z_p^m}\abs{f(\lambda,x)}\abs{dx},
\end{equation}
where $|dx|$ is the Haar measure on $\Z_p^m$ normalized so that
$\Z_p^m$ has measure $1$.
 \par
 A subanalytic constructible function on a subanalytic set $X$ is by definition
a $\Q$-linear combination of products of functions of the form
$v(h)$ and $|h'|$, where $h:X\to\Q_p^\times$ and $h':X\to\Q_p$ are
subanalytic functions. (For the notion of subanalytic functions
and subanalytic sets we refer to the section Terminology and
Notation below.)
 \par
We prove the following conjecture of Denef \cite{Denef1}:
\begin{theorem}\label{simple:1} The function $I$
is a subanalytic constructible function on $\Z_p^s$.
\end{theorem}
 \par In the case that the functions $f_i$ are polynomials, the map
$I$ has been studied by Igusa for $r=1$, by Lichtin for $r=2$, and
by Denef for arbitrary $r$ (see \cite{Igusa1, Igusa2, Igusa3},
\cite{Lichtin}, and \cite{Denef1}). In the more general case that
the $|f(\lambda,x)|$ in (\ref{integral}) is replaced by an
arbitrary subanalytic constructible function, the conclusion of
theorem~\ref{simple:1} still holds (see theorem~\ref{thm:basic}
below), where now $I$ takes the value zero whenever the integrated
function is not integrable.
 \par
The rationality of the analytic $p$-adic Serre - Poincar\'e series
was conjectured in \cite{Oest} and \cite{Serre} and proven by
Denef and van den Dries in \cite{DvdD}; the rationality can
immediately be obtained as a corollary of integration theorem
\ref{simple:1}. This is because it is well-known how to express
the Poincar\'e series as a $p$-adic integral (see \cite{Denef1},
section 1.6).
 \par
A second key result of the present paper is a cell decomposition
theorem for subanalytic sets and subanalytic functions
(theorem~\ref{Thm:Cell}), in perfect analogogy to the
semialgebraic cell decomposition theorem of \cite{Denef} and
\cite{Denef2}. Roughly speaking, $p$-adic cell decomposition
theorems describe the norm of given functions after partitioning
the domain of the functions in finitely many basic sets, called
cells. Cell decompositions are very useful to study parameterized
$p$-adic integrals (see below and \cite{Denef1}), and to prove the
rationality of Igusa's local zeta functions and of several
Poincar\'e series (see \cite{Denef}). Many of the applications of
cell decomposition (in for example ~\cite{Denef1} and
\cite{Denef3}) can, up to now, not be proven with other
techniques.
 \par
The proof of the analytic cell decomposition is based on several
results  by van den Dries, Haskell, and Macpherson \cite{vdDHM} on
the geometry of subanalytic $p$-adic sets and subanalytic
functions; we state some of these results in section \ref{proof
cell decomp}.
 \par
We also extensively use the theory of $p$-adic subanalytic sets,
developed by Denef and van den Dries in \cite{DvdD} in analogy to
the theory of real subanalytic sets; in particular, we use the
dimension theory of \cite{DvdD}. In section
\ref{section:classification} we apply cell decomposition to obtain
the following classification:
\begin{theorem}
\label{classification} Let $X\subset \Q_p^m$ and $Y\subset \Q_p^n$
be infinite subanalytic sets. Then there exists a subanalytic
bijection $X\to Y$ if and only if $\dim(X)=\dim(Y)$.
\end{theorem}
This classification of subanalytic sets is similar to the
classification of semialgebraic sets in \cite{C}. Note that in
particular there exists a semialgebraic bijection between $\Q_p$
and $\Q_p^\times$; this is the main result of \cite{CH}.
 \par
The theory of $p$-adic integration has also served as an inspiring
example for the theory of motivic integration and there are many
connections to it (see e.g.~\cite{DL} and \cite{DLinvent}).
 \par
Many of the results of \cite{DvdD} and \cite{vdDHM} are formulated
for $\Q_p$ and not for finite field extensions of $\Q_p$,
nevertheless, all results referred to in this paper, also hold for
finite field extensions of $\Q_p$ (see the remark in (3.31) of
\cite{DvdD}). All results of this paper also hold in finite field
extensions of $\Q_p$.
 \subsection*{Acknowledgment} I would like to thank Denef for
pointing out to me his conjecture on $p$-adic subanalytic
integrals and for having stimulating conversations on this and
related subjects. I thank van den Dries, Haskell, and Mourgues for
many interesting discussions. I also thank the referee for his
useful suggestions on the presentation of the paper.
 \subsection*{Terminology and Notation} Let $p$ denote a fixed prime number, $\Q_p$ the
field of $p$-adic numbers and $K$ a fixed finite field extension
of $\Q_p$ with valuation ring $R$. For $x\in K^\times$ let
$v(x)\in\Z$ denote the $p$-adic valuation of $x$ and
$|x|=q^{-v(x)}$ the $p$-adic norm, with $q$ the cardinality of the
residue class field. We write $P_n=\{y^n\mid y\in K^\times\}$ and
$\lP$ denotes $\{\mu x\mid x\in P_n\}$ for $\mu\in K$.
 \par
For $x=(x_1,\ldots,x_m)$ let $K\{ x\}$ be the ring of restricted
power series over $K$ in the variables $x$; it is the ring of
power series $\sum a_ix^i$ in $K[[x]]$ such that $|a_i|$ tends to
$0$ as $|i|\to\infty$. (Here, we use the multi-index notation
where $i=(i_1,\ldots,i_m)$, $|i|=i_1+\ldots+i_m$ and
$x^i=x_1^{i_1}\ldots x_m^{i_m}$.) For $x_0\in R^m$ and $f=\sum
a_ix^i$ in $K\{x\}$ the series $\sum a_ix_0^i$ converges to a
limit in $K$, thus, one can associate to $f$ a \emph{restricted
analytic function} given by
\[
f:K^m\to K:x\mapsto
 \left\{\begin{array}{ll} \sum_i a_i x^i & \mbox{ if
}x\in R^m,\\
0 & \mbox{ else.}
 \end{array}\right.
\]
We extend the notion of $D$-functions of \cite{DvdD} to our
setting\footnote{In \cite{DvdD}, $D$-functions are functions from
$R^m$ to $R$ for $m>0$.}:
\begin{definition} A $D$-function is a function $K^m\to K$ for some
$m\geq0$, obtained by repeated application of the following rules:
\begin{itemize}
 \item[(i)] for each $f\in K\{x_1,\ldots,x_m\}$, the associated restricted analytic
function $x\mapsto f(x)$ is a $D$-function;
 \item[(ii)] for each polynomial $f\in K[x_1,\ldots,x_m]$, the polynomial
 map $x\mapsto f(x)$ is a $D$-function;
 \item[(iii)] the function $x\mapsto x^{-1}$, where $0^{-1}=0$ by convention, is a
$D$-function;
 \item[(iv)] for each $D$-function $f$ in
$n$ variables and each $D$-functions $g_1,\ldots,g_n$ in $m$
variables, the function $f(g_1,\ldots,g_n)$ is a $D$-function.
\end{itemize}
\end{definition}
A \emph{(globally) subanalytic subset of} $K^m$  is a subset of
the form
\[X=\bigcup_{i=1}^r\bigcap_{j=1}^{s}X_{ij}\]
where each $X_{ij}$ is of the form $\{x\in K^m\mid f_{ij}(x)=0\}$
or $\{x\in K^m\mid f_{ij}(x)\in P_{n_{ij}}\}$, where the functions
$f_{ij}$ are $D$-functions and  $n_{ij}>0$. We call a function
$g:A\subset K^m\to K^n$ subanalytic if its graph is a subanalytic
set. We refer to \cite{DvdD}, \cite{Denef1}, and \cite{vdDHM} for
the theory of subanalytic $p$-adic geometry and to \cite{Lip} for
the theory of rigid subanalytic sets.
 \par
In section \ref{proof cell decomp} we will use the framework of
model theory. We let $\Ldan$ be the first order language
consisting of the symbols
\[
+,\ -,\ \cdot,\ ^{-1}, \{P_n\}_{n>0},
\]
together with an extra function symbol $f$ for each restricted
analytic function associated to restricted power series in
$\bigcup_m K\{x_1,\ldots,x_m\}$. We consider $K$ as an
$\Ldan$-structure using the natural interpretations of the symbols
of $\Ldan$.
 \par
We mention the following fundamental result in the theory of
subanalytic sets.
\begin{theorem}[\cite{DvdD}, Corollary (1.6)]\label{Gabrielov}
The collection of subanalytic sets is closed under taking
complements, finite unions, finite intersections, and images under
subanalytic maps.\footnote{I take the occasion to correct a small
error in \cite{vdDHM} with respect to quantifier elimination on
$\Q_p$. Namely, the division function $D$ in \cite{vdDHM} should
either be replaced by the field inverse $^{-1}$ or by the function
$D$ given by $D(x,y)=x/y$ if $0<|x|\leq|y|$ and $D(x,y)=0$
otherwise.}
 \end{theorem}
A \emph{semialgebraic subset of} $K^m$ is a subset of the same
form as $X$ above but with the $f_{ij}$ polynomials over $K$, and
a function is semialgebraic if its graph is a semialgebraic set.
It is well-known that also the collection of semialgebraic sets is
closed under taking complements, finite unions and intersections,
and images under semialgebraic maps (see \cite{Mac},
\cite{Denef2}).
 \section{Analytic cell decomposition}\label{proof
 cell decomp}
To state cell decomposition one needs basic sets called
(subanalytic) cells, which we define inductively. For $m,l>0$
write $\pi_m:K^{m+l}\to K^m$ for the linear projection on the
first $m$ variables and, for $A\subset K^{m+l}$ and $x\in
\pi_m(A)$, write $A_x$ for the fiber $\{t\in K^l\mid (x,t)\in
A\}$.
\begin{definition}\label{def:cell}
A cell $A\subset K$ is a (nonempty) set of the form
 \begin{equation}
\{t\in K\mid |\aaa|\sq_1 |t-\gamma|\sq_2 |\bbb|,\
  t-\gamma\in \lP\},
\end{equation}
with constants $n>0$, $\mu,\gamma\in K$, $\aaa,\bbb\in K^\times$,
and $\square_i$ either $<$ or no condition. If
$\mu=0$ we call $A$ a $0$-cell and we call $A$ a $1$-cell otherwise.\\
 A (subanalytic) cell $A\subset K^{m+1}$, $m\geq 1$, is a set  of the
form
 \begin{equation}\label{Eq:cell}
 \begin{array}{ll}
\{(x,t)\in K^{m+1}\mid
 &
 x\in D,\,  |\aaa(x)|\sq_1 |t-\gamma(
 x)|\sq_2 |\bbb(x)|,\\
 &
 t-\gamma(x)\in \lP\},
  \end{array}
\end{equation}
 with $(x,t)=(x_1,\ldots,
x_m,t)$, $n>0$, $\mu\in K$, $D=\pi_m(A)$ a cell, subanalytic
functions $\aaa,\bbb:K^m\to K^\times$, $\gamma:K^m\to K$, and
$\square_i$ either $<$ or no condition. We call $\gamma$ the
\emph{center} and $\mu P_n$ the \emph{coset} of the cell $A$. If
$D$ is a cell of type $(i_1,\ldots,i_m)$ with $i_j\in\{0,1\}$, we
call $A$ an $(i_1,\ldots,i_m,0)$-cell if $\mu=0$ and we call $A$
an $(i_1,\ldots,i_m,1)$-cell otherwise. If at each stage of this
inductive definition the functions $\aaa$, $\bbb$, and $\gamma$
are analytic on the respective projections  $\pi_i(A)$,
$i=1,\ldots,m-1$, we call $A$ an \emph{analytic cell}.
 \end{definition}
 \begin{remark}
 \item[(i)] An $(i_1,\ldots,i_m,0)$-cell $A\subset K^{m+1}$ is the graph of a
subanalytic function defined on $\pi_m(A)$, and, if $A$ is an
$(i_1,\ldots,i_m,1)$-cell, then $A_x$ is a nonempty open in $K$
for each $x\in \pi_m(A)$.
 \item[(ii)] An analytic cell is
a $K$-analytic manifold (in the obvious sense, see
e.g.~\cite{Bour}).
\end{remark}
\par
Theorem \ref{Thm:Cell} below is a subanalytic analogue of the
semialgebraic cell decomposition (see \cite{Denef} and
\cite{Denef2}); it is a perfect analogue of the reformulation
\cite[lemma 4]{C}. In \cite{Denef1}, an overview is given of
applications of $p$-adic cell decomposition, going from a
description of local singular series to counting profinite
$p$-groups (as in \cite{DuS}).
 \begin{theorem}[Analytic Cell Decomposition]\label{Thm:Cell}
Let $X\subset K^{m+1}$ be a subanalytic set, $m\geq0$, and
$f_j:X\to K$ subanalytic functions for $j=1,\ldots,r$. Then there
exists a finite partition of $X$ into cells $A$ with center
$\gamma:K^m\to K$ and coset $\mu P_n$ such that for each $(x,t)\in
A$
 \begin{equation}\label{eq:thm:cell}
 |f_j(x,t)|=
 |\ddd_j(x)|\, |(t-\gamma(x))^{a_j}\mu^{-a_j}|^\frac{1}{n},\qquad
 \mbox{ for each } j=1,\ldots,r,
 \end{equation}
with $(x,t)=(x_1,\ldots, x_m,t)$, integers $a_j$, and
$\ddd_j:K^m\to K$  subanalytic functions. If $\mu=0$, we use the
conventions $a_j=0$ and $0^0=1$. Moreover, the cells $A$ can be
taken to be analytic cells such that the $\delta_j$ are analytic
on $\pi_m(A)$.
 \end{theorem}
\begin{remark}
  \item[(i)]
Theorem \ref{Thm:Cell} can be seen as a $p$-adic analogue of the
preparation theorem \cite{LR} for real subanalytic functions, or
as an analogue of cell decomposition for real subanalytic sets
(see e.g.~\cite{vdD}).
 \item[(ii)] Some of the analytic analogues of applications in
\cite{Denef1} as well as some of the results of \cite{Denef3},
\cite{vdDHM} and \cite{DvdD} can be obtained immediately using
cell decomposition and the integration theorems of this paper, for
example: Cor.~1.8.2.~of \cite{Denef3} on local singular series,
Thm.~3.1 of \cite{Denef1}, the $p$-adic Lojasiewicz inequalities
(3.37), the subanalytic selection theorem (3.6), the
stratification theorem~(3.29), and Thm.~(3.2) of \cite{DvdD}.
However, note that the presented proof of Thm.~\ref{Thm:Cell}
relies on \cite{DvdD} and \cite{vdDHM}.
 \item[(iii)] Partial results towards subanalytic cell
decomposition have been obtained in \cite{Liu1},  \cite{Mourgues},
\cite{Mil} and \cite{Wilc}. In \cite{Mourgues} and \cite{Wilc}, a
partitioning of arbitrary subanalytic sets into cells is obtained,
but without the description of the norm of subanalytic functions
on these cells. However, the description of the norm of the
subanalytic functions in theorem \ref{Thm:Cell} is used in the
proof of the conjecture of Denef below.
\end{remark}
For the proof of theorem~\ref{Thm:Cell} we use techniques from
model theory, namely a compactness argument. (For general notions
of model theory we refer to \cite{Hodges}.)
 \par
Let $K_1$ be an $\Ldan$-elementary extension of $K$ and let $R_1$
be its valuation ring. In view of theorem~\ref{Gabrielov}, we can
call a set $X\subset K_1^m$ subanalytic if it is $\Ldan$-definable
(with parameters from $K_1$) and analogously for subanalytic
functions, cells, and so on. Expressions of the form $|x|<|y|$ for
$x,y\in K_1$ are abbreviations for the corresponding
$\Ldan$-formula's expressing $|x|<|y|$ for $x,y\in K$, as in lemma
2.1 of \cite{Denef2}\footnote{For example, if $K=\Q_p$ with
$p\not=2$, the property $|x|<|y|$ is equivalent to
$y^2+\frac{x^2}{p}\in P_2$.}. Cells in $K_1^m$ are defined just as
in $K^m$ by replacing everywhere $K$ by $K_1$ in the definition.
By a $D$-function $K_1^m\to K_1$ we mean a function given by an
$\Ldan$-term (with parameters from $K_1$) in $m$ variables.
Similarly, one can speak of semialgebraic subsets of $K_1^m$ (with
parameters from $K_1$)\footnote{This can be done using the
language of Macintyre, consisting of $+,-,\cdot,0,1,$ and the
collection of predicates $\{P_n\}$ for $n>0$.}.
 \par
We recall one of the main results of \cite{vdDHM}:
\begin{theorem}[theorem B of \cite{vdDHM}]
Each subanalytic subset of $K_1$ is semialgebraic.
\end{theorem}
The following two lemmas treat the one-dimensional part of
theorem~\ref{Thm:Cell}.
\begin{lemma}\label{lem:Dfunction}
Let $f:R_1\to K_1$ be a subanalytic function. Then there exists a
finite partition of $R_1$ into semialgebraic sets $A$ such that
for each $A$ there exist polynomials $p$ and $q$ such that
\[
|f(x)|=|p(x)/q(x)|^{1/e}, \mbox{ for each }x\in A,
\]
where $q$ has no zeros in $A$ and $e>0$ is an integer.
\end{lemma}
\begin{proof}
By theorem 3.6 of \cite{DenefBord}, there exists a finite
partition of $R_1$ into subanalytic sets $B$ such that
\[
|f(x)|=|g_B(x)/h_B(x)|^{1/e}, \mbox{ for each }x\in B,
\]
where $g_B$ and $h_B$ are $D$-functions, $h_B(x)\not=0$ on $B$ and
$e>0$. (In \cite{DenefBord} this is proven for subanalytic
functions $\Z_p^m\to \Z_p$ using quantifier elimination in an
elementary way; its proof can be repeated for our situation
$R_1\to K_1$ or otherwise one can instantiate parameters in the
result of \cite{DenefBord} to deduce this as a corollary.) By
Theorem B of \cite{vdDHM}, the sets $B$ are semialgebraic.
 \par
In \cite{vdDHM} it is proven that the norm of any $D$-function is
piecewise equal to the norm of a rational function, the pieces
being semialgebraic sets. More precisely, by proposition 4.1,
corollary 3.4 and lemma 2.10 of \cite{vdDHM}, there exists for
each function $g_B$ a finite partition of $R_1$ into semialgebraic
sets $C$ such that on each $C$
\[
|g_B(x)|=|g_{BC}(x)/h_{BC}(x)|, \mbox{ for each }x\in C,
\]
where $g_{BC}$ and $h_{BC}$ are polynomials over $K_1$ and
$h_{BC}(x)\not=0$ on $C$. The same holds for each function $h_B$.
Taking an appropriate partition using intersections of these sets
$C$ and $B$ the lemma follows.
\end{proof}
\begin{lemma}\label{cell decomp dim 1}
Let $X\subset R_1$ be a subanalytic set and $f:X\to K_1$ a
subanalytic function. Then there exists a finite partition $\PP$
of $X$ into cells, such that for each cell $A\in\PP$ with center
$\gamma\in K$ and coset $\mu P_n$
\[
|f(t)|=|\ddd|\, |(t-\gamma)^a\mu^{-a}|^\frac{1}{n}\mbox{ for each
} t\in A,
\]
with $\ddd\in K_1$ and $a$ an integer. We use the convention that
$a=0$ and $0^0=1$ when $\mu=0$.
\end{lemma}
\begin{proof} We extend $f$ to a function $R_1\to K_1$ by putting $f(x)=0$ if
$x\not\in X$. By theorem B of \cite{vdDHM}, the set $X$ is
semialgebraic. Apply lemma \ref{lem:Dfunction} to $f$ to obtain a
partition $\PP$. Intersecting each set in $\PP$ with $X$, we
obtain a partition $\PP'$ of $X$. Now apply the semialgebraic cell
decomposition (in the formulation of \cite[Lem.~4]{C}) to the sets
in $\PP'$ and the respective polynomials occurring in the
application of lemma \ref{lem:Dfunction}. If we refine the
obtained partition such that for each cell $A\subset C$ with coset
$\mu P_n$ the number $n$ is a multiple of $e$ (for the occurring
fractional powers $1/e$), then the lemma follows.
\end{proof}
 We will use the previous lemma and a model-theoretical compactness
argument to prove the following variant of theorem~\ref{Thm:Cell}.
\begin{theorem}
\label{thm:cell2} Let $K_1$ be an arbitrary $\Ldan$-elementary
extension of $K$ with valuation ring $R_1$. Let $X\subset
K_1^{m+1}$ be subanalytic and $f_j:X\to K_1$ subanalytic functions
for $j=1,\ldots,r$. Then there exists a finite partition of $X$
into subanalytic cells $A$ with center $\gamma:K_1^m\to K_1$ and
coset $\mu P_n$ such that for each $(x,t)\in A$
\[
 |f_j(x,t)|= |\ddd_j(x)|\,|(t-\gamma(x))^{a_j}\mu^{-a_j}|^\frac{1}{n},
 \]
with $(x,t)=(x_1,\ldots, x_m,t)$, integers $a_j$, and
$\ddd_j:K_1^m\to K_1$ subanalytic functions, $j=1,\ldots,r$. Here
we use the convention that $a_j=0$ and $0^0=1$ when $\mu=0$.
\end{theorem}
\begin{proof} The proof goes by induction on $m\geq 0$. It is enough to prove
the theorem for $r=1$ (the theorem then follows after a
straightforward further partitioning, see for example
\cite{Denef2}).
 \par
When $m=0$, the usual change of variables $t'=1/t$ reduces the
description of what happens outside $R_1$ to what happens on
$R_1$, and an application of lemma \ref{cell decomp dim 1} gives
the desired result.
 \par
Let $X\subset K_1^{m+1}$ and $f:X\to K_1$ be subanalytic, $m>0$.
We write $(x,t)=(x_1,\ldots,x_m,t)$ and know by the previous that
for each fixed $x\in K_1^m$ we can decompose the fiber $X_x$ and
the function $t\mapsto f(x,t)$ on this fiber. We will measure the
complexity of given decompositions on which $|f(x,\cdot)|$ has a
nice description and see that this must be uniformly bounded when
$x$ varies.
 \\
To do this, we define a countable set $\SSS=\{\lP\mid \mu\in K,\
n>0\}\times\Z\times \{<, \emptyset\}^2$ and
$\SSS'=(K_1^\times)^2\times K_1^2$. To each
$d=(\lP,a,\sq_1,\sq_2)$ in $\SSS$ and
$\xi=(\xi_1,\xi_2,\xi_3,\xi_4)\in \SSS'$ we associate a set
$Dom_{(d,\xi)}$ as follows
 \[
Dom_{(d,\xi)} = \{t\in K_1\mid |\xi_1|\sq_1 |t-\xi_3|\sq_2
|\xi_2|,\ t-\xi_3\in\lP\}
 \]
The set $Dom_{(d,\xi)}$ is either empty or a cell and is
independent of $\xi_4$ and $a$. For fixed $k>0$ and tuple
$d=(d_1,\ldots,d_k)\in \SSS^k$, let $\varphi_{(d,k)}(x,\xi)$ be an
$\Ldan$-formula in the free variables $x=(x_1,\ldots,x_m)$ and
$\xi=(\xi_1,\ldots,\xi_k)$, with
$\xi_i=(\xi_{i1},\xi_{i2},\xi_{i3},\xi_{i4})$, such that
$(x,\xi)\in K_1^{m+4k}$ satisfies $\varphi_{(d,k)}$  if and only
if the following are true:
\begin{itemize}
\item[(i)] $x\in\pi_m(X)$ and $\xi\in (\SSS')^k$,\\
\item[(ii)] the collection  of the sets $Dom_{(d_i,\xi_i)}$ for $i=1,\ldots,k$ forms
 a partition of the fiber $X_x=\{t\in K_1\mid (x,t)\in X\}$,\\
\item[(iii)]
$|\xi_{i4}|\,|(t-\xi_{i3})^{a_i}\mu_i^{-a_i}|^\frac{1}{n_i}=\abs{f(x,t)}$
for each  $t\in Dom_{(d_i,\xi_i)}$ and each $i=1,\ldots,k$.
 \end{itemize}
Now we define for each $k>0$ and $d \in \SSS^k$ the set
\[
B_d= \{x\in K_1^m\mid \exists \xi\quad \varphi_d(x,\xi)\}.
\]
Each set $B_d$ is subanalytic and the (countable) collection
$\{B_d\}_{k,d}$ covers $\pi_m(X)$, because each $x\in\pi_m(X)$ is
in some $B_d$ by the induction. Since $K_1$ is an arbitrary
elementary extension of $K$, finitely many sets of the form $B_d$
must already cover $\pi_m(X)$ by model-theoretical compactness.
Consequently, we can take subanalytic sets $D_1,\ldots,D_s$ such
that $\{D_i\}$ forms a partition of $\pi_m(X)$ and each $D_i$ is
contained in a set $B_d$ for some $k>0$ and $k$-tuple $d$. For
each $i=1,\ldots,s$, fix such a $d$ with $D_i\subset B_d$, and let
$\Gamma_i$ be the subanalytic set
\[
\Gamma_i=\{(x,\xi)\in D_i\times (\SSS')^k\mid \varphi_d(x,\xi)\}.
\]
Then $\pi_m(\Gamma_i)=D_i$ by construction ($\pi_m$ is the
projection on the $x$-coordinates). By theorem 3.6 \cite{DvdD} on
definable Skolem functions, there is a subanalytic function
$D_i\to K_1^{4k}$ associating to $x$ a tuple $\xi(x)\in (\SSS')^k$
such that $(x,\xi(x))\in\Gamma_i$ for each $x\in D_i$. The theorem
follows now by partitioning further with respect to the
$x$-variables and using the induction hypothesis.
 \end{proof}
 \begin{proof}[Proof of theorem \ref{Thm:Cell}]
We only have to show that we can partition $X$ using analytic
cells $A$ in  such a way that the functions $\delta_j$ are
analytic on $\pi_m(A)$. In \cite{DvdD} one proves that any
subanalytic function is piecewise analytic. Theorem \ref{Thm:Cell}
then follows from theorem \ref{Thm:Cell} by partitioning further
using this fact.
\end{proof}
 \section{Classification of Subanalytic Sets}\label{section:classification}
For $X\subset K^m$ subanalytic and nonempty, the dimension
$\dim(X)$ of $X$ is defined as the largest integer $n$ such that
there is a $K$-linear map $\pi:K^m\to K^n$ and a nonempty
$U\subset \pi(X)$, open in $K^n$ (for alternative definitions, see
\cite{DvdD}).
\begin{theorem}\label{prop:bijectie}
For any subanalytic set $X\subset K^m$ and subanalytic functions
$f_i:X\to K$, $i=1,\ldots,r$, there is a semialgebraic set $Y$, a
subanalytic bijection $F:X\to Y$ and there are semialgebraic maps
$g_i:Y\to K$ such that
\[|g_i( F(x))|=|f_i(x)|\qquad  \mbox{for each }x\in X.
\]
\end{theorem}
\begin{proof}
We will give a proof by induction on $m$. Suppose that $X\subset
K^{m+1}$ is subanalytic and that $f_i:X\to K$ are subanalytic
functions, $m\geq0$. Apply cell decomposition to $X$ and the
functions $f_i$ to obtain a finite partition $\PP$ of $X$. For
$A\in\PP$ and $(x,t)\in A$, suppose that
$|f_i(x,t)|=|\ddd_i(x)|\,|(t-\gamma(x))^{a_i}\mu^{-a_i}|^\frac{1}{n}$,
$i=1,\ldots,r$, and suppose that $A$ is a cell of the form
 \[
 \begin{array}{ll}
 \{(x,t)\in K^{m+1} \mid
 &
 x\in D,\, |\aaa(x)|\sq_1
 |t-\gamma(x)|\sq_2 |\bbb(x)|,\\
 &
 t-\gamma(x)\in\lP \},
 \end{array}
 \]
 like in (\ref{Eq:cell}).
 After the translation $(x,t)\mapsto (x,t-\gamma(x))$ we
may suppose that $\gamma$ is zero on $D$. Apply the induction
hypotheses to the sets $D$ and the subanalytic functions
$\aaa,\bbb$, and $\ddd_i$. Repeating this process for every $A\in
\PP$, and noting that there is a semialgebraic function $h:P_n\to
K$ such that $|h(t)|=|t|^{1/n}$, the proposition follows after
taking appropriate disjoint unions inside $K^m$ of the occurring
semialgebraic sets.
 \end{proof}
We prove the following generalization of Theorem
\ref{classification}.
\begin{theorem}
\label{thm:classification} Let $X\subset K^m$ and $Y\subset K^n$
be infinite subanalytic sets. Then there exists a subanalytic
bijection $X\to Y$ if and only if $\dim(X)=\dim(Y)$.
\end{theorem}
\begin{proof} By theorem \ref{prop:bijectie}
there are subanalytic bijections  $X\to X'$ and $Y\to Y'$ with
$X'$ and $Y'$ semialgebraic, but then there exists a semialgebraic
bijection $X'\to Y'$ if and only if $\dim(X')=\dim(Y')$  by
theorem 2 of \cite{C}. Since the dimension of a subanalytic set is
invariant under subanalytic bijections (see \cite{DvdD}), the
theorem follows.
 \end{proof}
\section{Parametrized Analytic Integrals}\label{section integrals}
We show that certain algebra's of functions from $\Q_p^m$ to the
rational numbers $\Q$  are closed under $p$-adic integration.
These functions are called subanalytic constructible functions and
they come up naturally when one calculates parametrized $p$-adic
integrals like (\ref{integral}).
 \par
For $x=(x_1,\ldots,x_m)$ an $m$-tuple of variables, we will write
$|dx|$ to denote the Haar measure on $K^m$, so  normalized that
$R^m$ has measure $1$.
\begin{definition}\label{basic algebra's}
For each subanalytic set $X$, we let $\OOL(X)$ be the $\Q$-algebra
generated by the functions $X\to\Q$ of the form $x\mapsto v(h(x))$
and $x\mapsto |h'(x)|$ where $h:X\to K^\times$ and $h':X\to K$ are
subanalytic functions. We call $f\in\OOL(X)$ a \emph{subanalytic
constructible function} on $X$.
 \par
To any function $f$ in $\OOL(K^{m+n})$, $m,n\geq 0$, we associate
a function $I_m(f):K^m\to \Q$ by putting
 \begin{equation}\label{I_l}
I_m(f)(\lambda)= \int\limits_{K^n}f(\lambda,x)|dx|
 \end{equation}
if the function $x\mapsto f(\lambda,x)$ is absolutely integrable
for all $\lambda\in K^m$, and  we put $I_m(f)(\lambda)=0$
otherwise.
\end{definition}
 \begin{theorem}[Basic Theorem on $p$-adic Analytic Integrals]\label{thm:basic}
For any function $f\in\OOL(K^{m+n})$, the function $I_m(f)$ is in
$\OOL(K^{m})$.
 \end{theorem}
 \begin{proof}
By induction and Fubini's theorem it is enough to prove that for a
function $f$ in $\OOL(K^{m+1})$ in the variables
$(\lambda_1,\ldots,\lambda_m,t)$ the function $I_m(f)$ is in
$\OOL(K^m)$. Suppose that $f$ is a $\Q$-linear combination of
products of functions $|f_i|$ and $v(g_j)$, $i=1,\ldots,r$,
$j=1,\ldots,s$ where $f_i$ and $g_j$ are subanalytic functions
$K^{m+1}\to K$ and $g_j(\lambda,t)\not=0$. Applying cell
decomposition to $K^{m+1}$ and the functions $f_i$ and $g_j$, we
obtain a partition $\cP$ of $K^{m+1}$ into cells such that
$I_m(f)(\lambda)$ is a sum of integrals over
$A_\lambda=\{t\mid(\lambda,t)\in A\}$ for each cell $A\in\cP$,
where the integrands on these pieces $A_\lambda$ have a very
simple form. More precisely, on each piece $A_\lambda$ the
integrand is a $\Q$-linear combination of functions of the form
\begin{equation}\label{int}
\ddd(\lambda) |(t-\gamma)^a\mu^{-a}|^\frac{1}{n}
v(t-\gamma(\lambda))^l
\end{equation}
where $A$ is a cell with center $\gamma:K^m\to K$ and coset $\mu
P_n$, and with integers $a$ and $0\leq l$, and a function $\ddd$
in $\OOL(K^m)$. We may suppose that $\ddd(\lambda)\not=0$ for some
$\lambda\in\pi_m(A)$. Regroup all such terms where the same
exponents $a$ and $l$ appear, possibly by replacing the functions
$\ddd(\lambda)$ by other functions in $\OOL(K^m)$ (their
respective sums). One checks that the integrability of such an
integrand then only depends on the integers $a,n$, and $l$
occurring in each of the terms as in (\ref{int}) and on the
symbols $\sq_i$ and $\mu$ occurring in the description of the cell
$A$. By consequence, we may suppose that the integrand is a single
term of the form like in (\ref{int}) and that this term is
absolutely integrable over $A$. It suffices to show that the
integral
\begin{equation}\label{intint}
\ddd(\lambda)\int\limits_{t\in A_\lambda}
|(t-\gamma(\lambda))^a\mu^{-a}|^\frac{1}{n}
v(t-\gamma(\lambda))^l|dt|
\end{equation}
is in $\OOL(K^m)$. Write $u=t-\gamma(\lambda)$; since $A$ is a
cell with center $\gamma$ and coset $\mu P_n$, the set $A$ is  of
the form
\[
A= \{(\lambda,u)\in K^{m+1}\mid   \lambda\in D, \
|\aaa(\lambda)|\sq_1 |u|\sq_2 |\bbb(\lambda)|,\
  u\in \mu P_n\},
\]
with $\sq_i$ either $<$ or no condition, $D$ a cell, and
$\aaa,\bbb:K^m\to K^\times$ subanalytic functions. Taking into
account that the integral (\ref{intint}) is, by supposition,
integrable, only a few possibilities can occur (with respect to
the integers $a,n,$ and $l$, the conditions $\sq_i$, and $\mu$).
If $\mu=0$, the set $A_\lambda$ is a point for each $\lambda\in
D$, thus the statement is clear. Suppose $\mu\not=0$. In case that
both $\sq_1$ and $\sq_2$ represent no condition, the integrand has
to be zero by the supposition of integrability, and the above
integral trivially is in $\OOL(K^m)$. We suppose from now on that
$\sq_1$ is $<$; the other cases can be treated similarly. The
integral (\ref{intint}) can be rewritten as
\[
\ddd(\lambda)\cdot\int\limits_{u\in A_\lambda}
|u^a\mu^{-a}|^\frac{1}{n} v(u)^l|du|
\]
\begin{eqnarray*}
& = &\ddd(\lambda)\sum_k (q^{-ak}|\mu^{-a}|)^\frac{1}{n}
k^{l}\cdot\mathrm{Measure}\{u\in
A_\lambda\mid v(u)=k\}\\
&= &\ \epsilon\ddd(\lambda) \sum_k (q^{-ak}|\mu^{-a}|)^\frac{1}{n}
k^{l}q^{-k}
\end{eqnarray*}
for $\epsilon=q^s\cdot\mathrm{Measure}\{u\in A_\lambda\mid
v(u)=s\}$ (where $s$ is any number such that $\emptyset\not=\{u\in
A_\lambda\mid v(u)=s\}$), and where the summation is over those
integers $k\equiv v(\mu) \bmod{n}$ satisfying
 \[
 |\aaa(\lambda)|<
q^{-k}\sq_2|\bbb(\lambda)|.
 \]
We may suppose that on $A$, the  residue classes
\[
v(\aaa(\lambda))\pmod{n}\quad\mbox{ and  }\quad
v(\bbb(\lambda))\pmod{n}
 \]
 are fixed (possibly after refining the partition $\PP$).
Then this sum is equal to a $\Q$-linear combination of products of
the functions $\ddd$, $|\aaa|$, $|\bbb|$, $v(\aaa)$ and $v(\bbb)$.
For example, if $a/n=-1$, $\sq_1$ and $\sq_2$ are necessarily $<$
and one obtains a polynomial in $v(\aaa)$ and $v(\bbb)$ of degree
$\leq l+1$, multiplied with $\delta$. For more examples of
calculations of sums of this kind, see \cite{Denef3}, proof of
lemma 3.2. Thus, the integral (\ref{intint}) is in $\OOL(K^m)$ as
was to be shown.
\end{proof}
As a corollary we will formulate another version of the basic
integration theorem, conjectured in \cite{Denef1} in the remark
following theorem~2.6.
\begin{definition}
A set $A\subset \Z^n\times K^m$ is called \emph{simple} if
\[\{(\lambda,x)\in K^{n+m}\mid (v(\lambda_1),\ldots,v(\lambda_n),x)\in
A\ \&\ \prod_{i=1\ldots,n}\ \lambda_i\not=0\}\]
 is a subanalytic set. A function $h:A\subset\Z^n\times \Q_p^m\to \Z$ is called \emph{simple} if
its graph is simple. For a simple set $X$ we let $\OOLp(X)$ be the
$\Q$-algebra generated by all simple functions on $X$ and all
functions of the form $q^h$ where $h$ is a simple function on $X$.
\par
For a function $f$ in $\OOLp(\Z^{k+l}\times K^{m+n})$, $k,l,m,n$
integers $\geq0$, we define $I_{k,m}(f):\Z^{k}\times K^{m}\to\Q$
as
\[
I_{k,m}(f)(z,\lambda)=\sum_{z'\in
\Z^{l}}\int\limits_{K^n}f(z,z',\lambda,x)|dx|
\]
if the function $(z',x)\mapsto f(z,z',\lambda,x)$ is absolutely
integrable for all $(z,\lambda)\in \Z^k\times K^m$ with respect to
the Haar measure on $K^n$ and the discrete measure on $\Z^l$, and
we define $I_{k,m}(f)(z,\lambda)=0$ otherwise.
\end{definition}
\begin{theorem}\label{Thm:conjecture}
For each $f$ in $\OOLp(\Z^{k+l}\times K^{m+n})$, the function
$I_{k,m}(f)$ is in $\OOLp(\Z^k\times K^l)$.
\end{theorem}
\begin{proof}
It is enough to prove that for a function $f$ in $\OOLp(\Z^k\times
K^{m})$ in the variables $(z_1,\ldots,z_k,x_1,\ldots,x_m)$ the
function obtained by eliminating $x_m$ by integration, resp.
eliminating $z_k$ by summation, is in the respective algebra
$\OOLp$.
 \par
We first focus on integration with respect to $x_m$. To
$f:\Z^k\times K^{m}\to\Q$ we can associate a function
$g:K^{k+m}\to \Q$ by replacing the variables $z$ running over
$\Z^k$ by variables $\lambda$ running over $K^k$ in such a way
that $g(\lambda,x)=f(v(\lambda_1),\ldots,v(\lambda_k),x)$ for each
$\lambda\in (K^\times)^k$ and $g(\lambda,x)=0$ if one of the
$\lambda_i$ is zero. By the definitions it is immediate that $g$
is in $\OOL(K^{k+m})$ and the integral of $f$ with respect to
$x_m$ corresponds to the integral of the function $g$ with respect
to $x_m$. If we eliminate $x_m$ by integration from $g$, then we
get the function $I_{k+m-1}(g)$ which is in $\OOL(K^{k+m-1})$ by
theorem \ref{thm:basic}. This function only depends on
$(v(\lambda_1),\ldots,v(\lambda_k),x_1,\ldots,x_{m-1})$ and thus
corresponds to a function in $\OOLp(K^{k+m-1})$ as one can check
(for example by using cell decomposition again).
 \par
If we want to eliminate $z_k$ by summation, we associate to $f$
the subanalytic constructible function $g':K^{k+m}\to \Q$
determined by
\[
g'(\lambda,x)=|\lambda_k|^{-1}\frac{p}{p-1}f(v(\lambda_1),\ldots,v(\lambda_k),x)
\]
if $\prod_{i=1\ldots,n}\ \lambda_i\not=0$ and $g'(\lambda,x)=0$ if
$\prod_{i=1\ldots,n}\ \lambda_i=0$. Integrating with respect to
$\lambda_k$ then corresponds to summing over $z_k$, and the same
argument as above can be applied to complete the proof.
 \end{proof}
\bibliographystyle{amsplain}
\bibliography{anbib}

\providecommand{\bysame}{\leavevmode\hbox to3em{\hrulefill}\thinspace}
\providecommand{\MR}{\relax\ifhmode\unskip\space\fi MR }
\providecommand{\MRhref}[2]{%
  \href{http://www.ams.org/mathscinet-getitem?mr=#1}{#2}
}
\providecommand{\href}[2]{#2}
\begin{thebibliography}{10}

\bibitem{Bour}
N.~Bourbaki, \emph{Vari\'et\'es diff\'erentielles et analytiques. {F}ascicule
  de r\'esultats}, Hermann, Paris, 1967, (French).

\bibitem{C}
R.~Cluckers, \emph{Classification of semialgebraic $p$-adic sets up to
  semialgebraic bijection}, Journal f{\"u}r die reine und angewandte Mathematik
  \textbf{540} (2001), 105--114.

\bibitem{CH}
R.~Cluckers and D.~Haskell, \emph{{G}rothendieck rings of $\mathbb{Z}$-valued
  fields}, Bulletin of Symbolic Logic \textbf{7} (2001), no.~2, 262--269.

\bibitem{Denef}
J.~Denef, \emph{The rationality of the {P}oincar\'e series associated to the
  $p$-adic points on a variety}, Inventiones Mathematicae \textbf{77} (1984),
  1--23.

\bibitem{Denef3}
\bysame, \emph{On the evaluation of certain $p$-adic integrals}, Th\'eorie des
  nombres, S\'emin. Delange-Pisot-Poitou 1983--84, vol.~59, 1985, pp.~25--47.

\bibitem{Denef2}
\bysame, \emph{$p$-adic semialgebraic sets and cell decomposition}, Journal
  f{\"u}r die reine und angewandte Mathematik \textbf{369} (1986), 154--166.

\bibitem{DenefBord}
\bysame, \emph{Multiplicity of the poles of the poincar\'e series of a $p$-adic
  subanalytic set}, S\'em. Th. Nombres Bordeaux \textbf{43} (1987-1988), 1--8.

\bibitem{Denef1}
\bysame, \emph{Arithmetic and geometric applications of quantifier elimination
  for valued fields}, MSRI Publications, vol.~39, pp.~173--198, Cambridge
  University Press, 2000.

\bibitem{DvdD}
J.~Denef and {L. van den} Dries, \emph{$p$-adic and real subanalytic sets},
  Annals of Mathematics \textbf{128} (1988), no.~1, 79--138.

\bibitem{DLinvent}
J.~Denef and F.~Loeser, \emph{Germs of arcs on singular algebraic varieties and
  motivic integration}, Inventiones Mathematicae \textbf{135} (1999), 201--232.

\bibitem{DL}
\bysame, \emph{Definable sets, motives and $p$-adic integrals}, Journal of the
  American Mathematical Society \textbf{14} (2001), no.~2, 429--469.

\bibitem{vdD}
{L. van den} Dries, \emph{Tame topology and o-minimal structures}, Lecture note
  series, vol. 248, Cambridge University Press, 1998.

\bibitem{vdDHM}
{L. van den} Dries, D.~Haskell, and D.~Macpherson, \emph{One-dimensional
  $p$-adic subanalytic sets}, Journal of the London Mathematical Society
  \textbf{59} (1999), no.~1, 1--20.

\bibitem{DuS}
M.P.F. du~Sautoy, \emph{Finitely generated groups, $p$-adic analytic groups and
  {P}oincar\'e series}, Annals of Mathematics \textbf{137} (1993), no.~3,
  639--670.

\bibitem{Hodges}
W.~Hodges, \emph{Model theory}, Encyclopedia of Mathematics and Its
  Applications, vol.~42, Cambridge University Press, 1993.

\bibitem{Igusa1}
J.~Igusa, \emph{Complex powers and asymptotic expansions {I}}, Journal f{\"u}r
  die reine und angewandte Mathematik \textbf{268} (1974), 110--130.

\bibitem{Igusa2}
\bysame, \emph{Complex powers and asymptotic expansions {II}}, Journal f{\"u}r
  die reine und angewandte Mathematik \textbf{278} (1975), 307--321.

\bibitem{Igusa3}
\bysame, \emph{Lectures on forms of higher degree (notes by {S}. {R}aghavan)},
  Lectures on mathematics and physics, Tata institute of fundamental research,
  vol.~59, Springer-Verlag, 1978.

\bibitem{Lichtin}
B.~Lichtin, \emph{On a question of {I}gusa: towards a theory of several
  variable asymptotic expansions {I}}, Compositio Mathematica \textbf{120}
  (2000), no.~1, 25--82.

\bibitem{LR}
J.-M. Lion and J.-P. Rolin, \emph{Int{\'e}gration des fonctions
  sous-analytiques et volumes des sous-ensembles sous-analytiques. (integration
  of subanalytic functions and volumes of subanalytic subspaces)}, Ann. {I}nst.
  {F}ourier \textbf{48} (1998), no.~3, 755--767, (French).

\bibitem{Lip}
L.~Lipshitz and Z.~Robinson, \emph{Rings of separated power series and
  quasi-affinoid geometry}, Paris: Soci\'et\'e Math\'ematique de France, vol.
  264, Ast\'erisque, 2000.

\bibitem{Liu1}
Nianzheng Liu, \emph{Analytic cell decomposition and the closure of $p$-adic
  semianalytic sets}, Journal of Symbolic Logic \textbf{62} (1997), no.~1,
  285--303.

\bibitem{Mac}
A.~Macintyre, \emph{On definable subsets of $p$-adic fields}, Journal of
  Symbolic Logic \textbf{41} (1976), 605--610.

\bibitem{Mourgues}
M.-H. Mourgues, \emph{Corps p-minimaux avec fonctions de {S}kolem
  d{\'e}finissables}, S{\'e}minaire de structures alg{\'e}briques
  ordonn{\'e}es, 1999--2000, pr\'epublication de l'\'equipe de logique
  math\'ematique de Paris 7, pp.~1--8.

\bibitem{Mil}
A.~Mylnikov, \emph{$p$-adic subanalytic preparation and cell decomposition
  theorems}, Ph.D. thesis, Purdue University, 1999.

\bibitem{Oest}
J.~Oesterl\'e, \emph{R\'eduction modulo $p^n$ des sous-ensembles analytiques
  ferm\'es de $\mathbb{Z}_p^n$}, Inventiones Mathematicae \textbf{66} (1982),
  no.~2, 325--341.

\bibitem{Serre}
J.-P. Serre, \emph{Quelques applications du th\'eor\`eme de densit\'e de
  {C}hebotarev}, Publ. Math. Inst. Hautes \'Etudes Sci. \textbf{323-401}
  (1981).

\bibitem{Wilc}
S.~Wilcox, \emph{Topics in the model theory of $p$-adic numbers}, Ph.D. thesis,
  University of Oxford, (unfinished).

\end{thebibliography}
\end{document}